\documentclass[11pt,notitlepage,twoside,a4paper]{amsart}
 \usepackage{amsfonts}

\usepackage{amsmath,amssymb,enumerate}

\usepackage{epsfig,fancyhdr,color}

\usepackage{amssymb}
\usepackage{amsmath,amsthm}
\usepackage{latexsym}
\usepackage{amscd}
\usepackage{psfrag}
\usepackage{graphicx}
\usepackage[latin1]{inputenc}
\usepackage[mathcal]{eucal}

\renewcommand{\textsc}{\textcolor{red}}

\newtheorem{theorem}{\rm\bf Theorem}[section]
\newtheorem{proposition}[theorem]{\rm\bf Proposition}
\newtheorem{lemma}[theorem]{\rm\bf Lemma}
\newtheorem{corollary}[theorem]{\rm\bf Corollary}
\newtheorem*{theorem 1}{\rm\bf Proposition 1}
\newtheorem*{theorem 2}{\rm\bf Proposition 2}

\theoremstyle{definition}

\theoremstyle{remark}
\newtheorem{remark}[theorem]{\rm\bf Remark}

\newtheorem{question}[theorem]{\rm\bf Question}

\def\interieur#1{\mathord{\mathop{\kern 0pt #1}\limits^\circ}}

\title[]{Shortening all the simple closed geodesics on surfaces with boundary}

\author{Athanase Papadopoulos}
\address{Athanase Papadopoulos, Max-Plank-Institut f\"ur Mathematik, Vivatsgasse 7, 53111 Bonn, Germany, and : Institut de Recherche Math{\'e}matique Avanc\'ee,
Universit{\'e} de Strasbourg and CNRS,
7 rue Ren\'e Descartes,
 67084 Strasbourg Cedex, France} \email{papadopoulos@math.u-strasbg.fr}
\date{\today}

\author{Guillaume Th\'eret}
\address{Guillaume Th\'eret, Max-Plank-Institut f\"ur Mathematik, Vivatsgasse 7, 53111 Bonn, Germany}
\email{theret@mpim-bonn.mpg.de}

\begin{document}

\begin{abstract}  
We give a proof of an unpublished result of Thurston showing that given any hyperbolic metric on a  surface of finite type with nonempty boundary, there exists another hyperbolic metric on the same surface for which the lengths of all simple closed geodesics are shorter.  (This is not possible for surfaces of finite type with empty boundary.) Furthermore, we show that we can do the shortening in such a way that it is bounded below by a positive constant. This improves a recent result obtained by Parlier in \cite{Parlier}. 
We include this result in a discussion of the weak metric theory of the Teichm\"uller space of  surfaces with nonempty boundary. \bigskip

\medskip

 
\end{abstract}
\maketitle

\section{Introduction}
\label{intro}

Let $S$ be a connected orientable surface of finite topological type and of negative Euler characteristic. All the hyperbolic structures that we shall consider on $S$ are metrically complete, and have finite area with totally geodesic boundary. 
Unless explicitly specified, we shall assume that the boundary $\partial S$ of $S$ is nonempty and that the boundary components are all closed geodesics.

Let $\mathcal{T}(S)$ denote the Teichm\"uller space of $S$, that is, 
the space of hyperbolic structures on $S$ of the type specified above up to homeomorphisms homotopic to the identity. In this paper, all homotopies of a surface fix setwise (and not necessarily pointwise) the boundary components.

Let  $\mathcal{C}=\mathcal{C}(S)$ be the set of simple closed geodesics in $S$, boundary components included. 
This set is defined relative to a hyperbolic structure which is understood, 
and it is known that there exists a natural correspondence between any two such sets relative to different underlying hyperbolic structures.

A {\it weak metric} on a set is a structure that satisfies all the axioms of a metric except the symmetry axiom.

We consider the following function on $\mathcal{T}(S)\times\mathcal{T}(S)$:
\begin{equation}\label{eq:k}
k(X,Y)=\log\sup_{\gamma\in\mathcal{C}}\frac{l_{Y}(\gamma)}{l_{X}(\gamma)}.
\end{equation}
A result of Thurston (obtained by combining Theorem 3.1 and Proposition 3.5 of \cite{Thurston1986}) says that in the case where $\partial S$ is empty, 
the function $k$ defines a weak metric on the Teichm\"uller space of $S$.

It is easy to see that for surfaces with nonempty boundary, the function $k$ is not a weak metric. Indeed, as already remarked in \cite{Liu-Papad-Su-Th2009}, it suffices to take $S$ to be a pair of pants (a sphere with three boundary components) and 
 $X$ and $Y$ two hyperbolic metrics on $S$ such that the lengths of the three boundary components for the metric 
$X$ are all strictly smaller than the corresponding lengths for the metric $Y$. It is clear that in this case we have $k(X,Y)<0$. We asked in the same paper whether this example of the pair of pants can be generalized to any surface with boundary, that is, if for any hyperbolic metric on a surface with  nonempty boundary, there exists another hyperbolic metric
on the same surface for which the lengths of all the simple closed geodesics is strictly decreased by a uniformly bounded amount.  Theorem \ref{th:shrink} that we prove below answers positively this question. In particular,
 the function $k$ is not a weak metric, for any surface $S$ with nonempty boundary.\\

We shall call a {\it simple geodesic arc} in $S$ a geodesic segment which is properly embedded in that surface, 
that is, the arc has no self-intersection, the interior of the arc is in the interior of $S$ and the endpoints of the arc are on $\partial S$.
Let $\mathcal{B}=\mathcal{B}(S)$ be the union of the set of geodesic boundary components of $S$ 
with the set of simple geodesic arcs that are perpendicular to the boundary. 
(The same remark as for the set $\mathcal{C}$ holds, namely, the set $\mathcal{B}$ is defined relative to some hyperbolic structure, 
but there exists a natural correspondence between two such sets relative to different hyperbolic structures.)

In contrast with the function $k$ defined in (\ref{eq:k}), we proved in \cite{Liu-Papad-Su-Th2009} that the function
$$
K(X,Y)=\log\sup_{\gamma\in\mathcal{B}}\frac{l_{Y}(\gamma)}{l_{X}(\gamma)}.
$$
is a weak metric on Teichm\"uller space.

\noindent {\it Acknowledgment.} We are grateful to the referee of this paper for several useful remarks and corrections.

\section{Shrinking all simple closed geodesics}

\begin{figure}[!hbp]
\centering
\includegraphics[width=.6\linewidth]{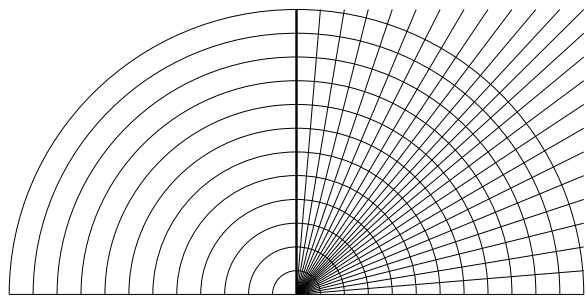}
\caption{\small{In the upper-half plane model of the hyperbolic plane, 
the foliation by Euclidean half-circles is the foliation $\mathcal{G}$ by geodesic lines, and the orthogonal foliation in the quarter plane to the right is the orthogonal foliation $\mathcal{E}$ 
by equidistant lines  to the vertical geodesic line. (Only the part of $\mathcal{E}$ in the right-quarter plane is drawn.)}}
\label{equid2}
\end{figure}

Consider a foliation $\mathcal{E}$ of the hyperbolic plane $\mathbb{H}^2$ by the set of curves that are equidistant from a given geodesic, and
consider the foliation $\mathcal{G}$ of $\mathbb{H}^2$ by the curves that are orthogonal to the leaves of $\mathcal{E}$ (Figure \ref{equid2}). 
The leaves of $\mathcal{G}$ are geodesics. 
We start with the following:

\begin{lemma}[Projection along equidistant curves] \label{lem:proj}
The projection map from $\mathbb{H}^2$ to some leaf of $\mathcal{G}$ along the leaves of $\mathcal{E}$ is distance non-increasing. Furthermore, the distance between any two points in $\mathbb{H}^2$ is equal to the distance between their projections if and only if the two points are on some leaf of $\mathcal{G}$.
\end{lemma}

\begin{proof}
Let $\gamma$ be the geodesic line in $\mathbb{H}^2$ of which $\mathcal{E}$ is the set of equidistant curves. 
The curves that are orthogonal to the curves of the family $\mathcal{E}$ are precisely the geodesic lines that are orthogonal to $\gamma$.

Let $x,y$ be two points in the hyperbolic plane $\mathbb{H}^{2}$.
If these points lie on the same leaf of $\mathcal{E}$, their projection is a point, and the result follows in this case.
Thus we can assume that the points $x$ and $y$ lie on distinct leaves of the foliation $\mathcal{E}$.
Consider the geodesic segment, $\sigma$, joining $x$ to $y$.
It may happen that the geodesic segment $\sigma$ touches some leaf of $\mathcal{E}$, but such a point of tangency is then unique. By dividing $\sigma$ into two geodesic segments that meet at that tangency point, we can assume without loss of generality, in the proof that follows, that $\sigma$ is transverse to the leaves of $\mathcal{E}$.
The goal is to compare the length of $\sigma$ with the length of any geodesic arc that is perpendicular to $\mathcal{E}$ 
and whose endpoints lie on the same equidistant curves as the endpoints of $\sigma$.
If the segment $\sigma$ is itself contained in a leaf of $\mathcal{G}$, 
then the projection of $\sigma$ keeps the length of $\sigma$ constant.
Thus, we can assume that $\sigma$ is not contained in a leaf of $\mathcal{G}$.

Up to dividing $\sigma$ into two geodesic segments, we can assume that the interior of $\sigma$ lies in a single
component of $\mathbb{H}^{2}\setminus\gamma$.
Furthermore, since the geodesic arcs on which we project $\sigma$ have all the same length, 
we can assume that the geodesic arc, $k$, on which we project $\sigma$ has a unique endpoint in common with $\sigma$.

There are two possibilities for choosing the arc $k$, which correspond
to the two possibilities for the common endpoint between $k$ and $\sigma$.
Let us specify a choice for this common endpoint and let us call it $A$.
Consider the two leaves of $\mathcal{E}$ passing through $x$ and $y$.
Since the segment $\sigma$ lies in a single component of $\mathbb{H}^{2}\setminus\gamma$, 
one of these two leaves is farther from $\gamma$ than the other.
Let us choose the arc $k$ so that the common point $A$ lies on this farthest leaf.

\begin{figure}[!hbp]
\psfrag{A}{$A$}
\psfrag{B}{$B$}
\psfrag{B'}{$B'$}
\psfrag{C}{$C$}
\psfrag{O}{$O$}
\psfrag{x}{}
\psfrag{y}{}
\psfrag{z}{}
\psfrag{h}{$\sigma$}
\psfrag{a}{$\gamma$}
\psfrag{g}{$k$}
\centering
\includegraphics[width=1.1\linewidth]{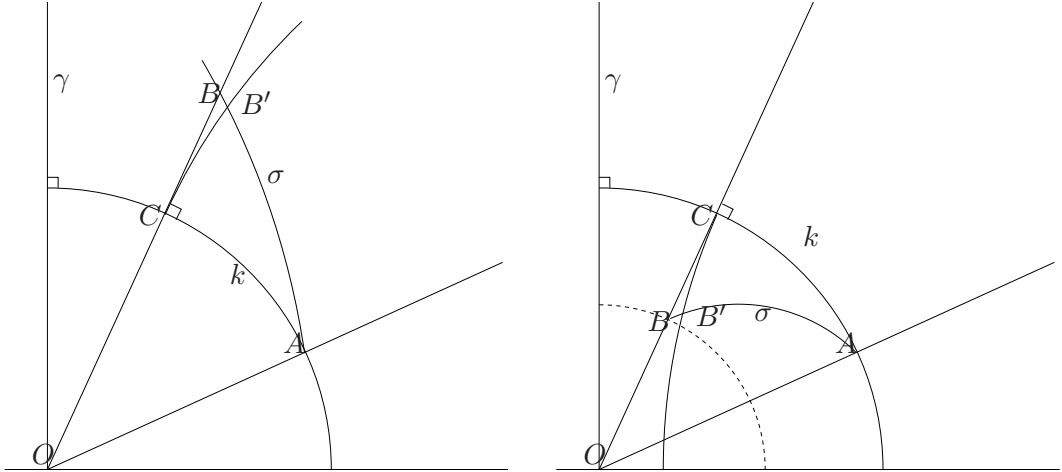}
\caption{\small{In these two figures we have represented a geodesic segment $\sigma$ which is transverse
to the foliation $\mathcal{E}$ whose leaves are equidistant curves from the geodesic $\gamma$. 
There are two natural candidates for the geodesic segment $k$ onto which one can project $\sigma$.
We consider the one for which the geodesic segment perpendicular to $k=[AC]$ through $C$ cuts $\sigma$ in an interior point.
In the left-hand side picture, $k$ lies below $\sigma$ whereas it lies above $\sigma$ in the right-hand side picture.}}
\label{equid1}
\end{figure}

We are led to consider the curvilinear triangle $ABC$ having two geodesic edges, 
namely $[AB]=\sigma$ and the geodesic segment $[AC]=k$ on which we project $\sigma$, 
and whose third edge $[BC]$ is the arc of an equidistant curve that connects
the endpoint $B$ of $\sigma$ to the endpoint $C$ of $k$ (see Figure \ref{equid1}).
Note that the angle $\widehat{BCA}$ is a right angle.  

Consider the geodesic passing through $C$ perpendicularly to $k$.
By convexity and thanks to our choice for $k$ (see Figure \ref{equid1}), 
this geodesic intersects the segment $\sigma$ in an interior point $B'$.
Hence,  
$$
AB\geq AB'.
$$
Now consider the geodesic triangle $AB'C$.
Since the angle $\widehat{B'CA}$ is a right angle,
we have, by hyperbolic trigonometry,
$$
\cosh(AB')=\cosh(B'C)\cosh(AC).
$$ 
Hence, 
$$
AB'>AC,
$$
since $B'C>0$ by assumption.
Thus, we have
$$
AB>AC,
$$
that is, the length of $\sigma$ is strictly greater than the length of $k$.
All the cases have been dealt with.
This concludes the proof.
\end{proof}

We recall a few facts about Nielsen extensions of hyperbolic surfaces with boundary.
Let $X$ be a hyperbolic structure on $S$.
With the above requirements on $S$, the {\it Nielsen extension} $\hat{X}$ of $X$ can be defined as the unique complete hyperbolic surface without boundary which contains $X$ and which retracts on $X$.
Another description of $\hat{X}$  is that this surface is obtained from $X$ by gluing a {\it funnel}, that is, a semi-infinite cylinder with one geodesic boundary, along each boundary components of $X$.
Note that the isometry type of each semi-infinite cylinder we glue is completely determined by the length of its unique boundary component, 
and that the hyperbolic structure $\hat{X}$ does not depend upon the way these cylinders are glued to $\partial S$ (that is, the twist parameters have no contribution).
Note also that the hyperbolic surface $\hat{X}$ has infinite area.

Let us remark that the Nielsen extension $\hat{X}$, 
although it is a naturally defined complete hyperbolic structure on a surface homeomorphic to the interior to $S$, 
is not isometric to the unique (Poincar\'e) complete hyperbolic structure on the interior of $S$ 
that is in the conformal class of the restriction of the metric $X$ to that interior. 

The {\it convex core} of a complete infinite-area hyperbolic structure on a surface of finite type is the hyperbolic surface with boundary 
obtained by cutting out each infinite half-cylinder along the unique geodesic on which it naturally retracts. 
The convex core of the Nielsen extension $\hat{X}$ of $X$ is the hyperbolic surface $X$ that we started with.

At the level of the universal coverings, we have the following picture:
The universal covering of the hyperbolic surface $X$ with boundary is a subset of the hyperbolic plane
bounded by the preimage of the boundary $\partial S$.
This preimage consists in infinitely many disjoint geodesic lines.
(If one identifies the hyperbolic plane with the unit disk, 
the limit set of the corresponding Fuchsian group is a Cantor set of the unit circle.)
The universal covering of the Nielsen extension $\hat{X}$ of $X$ is  
the hyperbolic plane $\mathbb{H}^2$, and it naturally contains the universal covering of $X$.
The infinite half-cylinders in $\hat{X}$ lift to the closed half-planes in the complement of the universal covering of $X$.\\

\begin{figure}[!hbp]
\psfrag{g1}{$\widetilde{\gamma_1}$}
\psfrag{g2}{$\widetilde{\gamma_2}$}
\psfrag{S}{$S_{\epsilon}$}
\centering
\includegraphics[width=.4\linewidth]{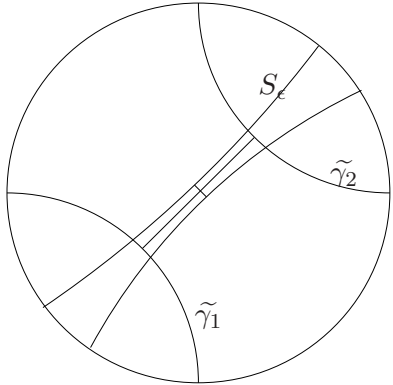}
\caption{\small{The segment $\widetilde{\alpha}$ is the geodesic segment joining perpendicularly the two hyperparallel geodesics $\widetilde{\gamma_1}$ and $\widetilde{\gamma_2}$. The 
 $\epsilon$-strip $S_{\epsilon}$ is bounded by two other hyperparallel geodesics drawn.}}
\label{equid3}
\end{figure}

Consider two hyperparallel geodesic lines in $\mathbb{H}^2$ and 
let $\widetilde{\alpha}$ be their common perpendicular geodesic segment.
Let $\epsilon$ be a positive number.
An {\it $\epsilon$-strip} $S_{\epsilon}$ around $\widetilde{\alpha}$ is a strip containing $\widetilde{\alpha}$ and bounded by two hyperparallel geodesics
whose common perpendicular has length $\epsilon$ (see Figure \ref{equid3}).
The {\it core} of the $\epsilon$-strip $S_{\epsilon}$ is this common geodesic segment joining the boundary components of $S_{\epsilon}$ perpendicularly and which is perpendicular to $\widetilde{\alpha}$ at their common midpoint.
Note that the core, $c_{\epsilon}$ of $S_{\epsilon}$ has length $\epsilon$.
We shall equip an $\epsilon$-strip with the foliation by arcs that are equidistant from the core. 
This foliation induces an isometric correspondence between the boundary geodesics of the $\epsilon$-strip, which we shall refer to as the {\it canonical} isometry between these geodesics.

Let $\alpha$ be a simple geodesic arc joining perpendicularly a boundary component $\gamma_{1}$ of the hyperbolic surface $X$ to a boundary component $\gamma_{2}$.
(We may have $\gamma_{1}=\gamma_{2}$.)
Consider lifts $\widetilde{\gamma}_{1}$, $\widetilde{\gamma}_{2}$ of $\gamma_{1}$, $\gamma_{2}$ to the universal covering.
The lifts $\widetilde{\gamma}_{1}$, $\widetilde{\gamma}_{2}$ are hyperparallel geodesic lines and  
there is a unique lift $\widetilde{\alpha}$ of $\alpha$ that joins them perpendicularly.
For any small enough $\epsilon>0$, 
the $\epsilon$-strip around $\widetilde{\alpha}$ projects to an embedded strip containing $\alpha$ in the 
Nielsen extension $\hat{X}$ of $X$.
We call such a strip an $\epsilon$-strip in $\hat{X}$.

We now define a construction that we call {\it peeling an $\epsilon$-strip from a hyperbolic surface with boundary}. 
Start from a hyperbolic structure $X$ on the surface with boundary $S$ and consider $X$ as embedded in its Nielsen extension $\hat{X}$.
Consider an $\epsilon$-strip, $B$, in $\hat{X}$.
Consider the hyperbolic structure $\hat{Y}_{B}$ on $\hat{X}$ obtained by cutting out the strip $B$ from $\hat{X}$ and 
by gluing back the geodesic sides of the closure of $\hat{X}\setminus B$ by the canonical isometry that identifies the endpoints of the core of $B$.
Another way of obtaining $\hat{Y}_{B}$ is by collapsing the strip $B$ along the leaves of the foliation of this strip by equidistant arcs defined above.
Let $\hat{f}_{B}\,:\,\hat{X}\to\hat{Y}_{B}$ be the collapsing map. 
Let $Y_{B}$ be the hyperbolic structure on $S$ obtained by restricting the hyperbolic structure $\hat{Y}_{B}$ 
to its convex core.  

The image, $\hat{\alpha}_{B}$ of the strip $B$ by $\hat{f}_{B}$ is an infinite geodesic. 

\begin{proposition}
\label{le:decrease}
The map $\hat{f}_{B}\,:\,\hat{X}\to\hat{Y}_{B}$ is $1$-Lipschitz and it is homotopic to the identity map of $S$.
More precisely, $\hat{f}_{B}$ is length-preserving in the complement of $B$ and it strictly decreases, by a uniform amount, 
distances between points that are ``separated by" $B$ and contained in $X$, that is, points that
can be joined by a geodesic of shortest length which intersects the strip $B$ transversely and is contained in 
$X$.
\end{proposition}

\begin{proof}
The assertion regarding the homotopy between $\hat{f}_{B}$ and the identity map is clear. (Note that the natural retractions from $\hat{X}$ to $X$ and from  $\hat{X}$ to $\hat{Y}_{B}$ allow us to talk about homotopies on $S$.)
Also, it is clear that the map $\hat{f}_{B}$ restricted to $\hat{X}\setminus B$ is length-preserving. Thus, it remains to show that the map strictly decreases the distances between 
points separated by $B$ by a uniform amount.

Let $x,y\in X\subset\hat{X}$ be two points separated by $B$ and let $[xy]$ be a shortest geodesic segment joining $x$ to $y$.
The length $xy$ of this segment is equal to $d_{\hat{X}}(x,y)$.
By the assumptions made on $x$ and $y$, the intersection $[xy]\cap B$ has only one component, and we denote it by $[x'y']$.
The image of $[xy]$ by $\hat{f}_{B}$ is a piecewise geodesic curve, namely,
$[\hat{f}_{B}(x)\hat{f}_{B}(x')]\cup[\hat{f}_{B}(x')\hat{f}_{B}(y')]\cup[\hat{f}_{B}(y')\hat{f}_{B}(y)]$.
Hence, 
$$
d_{\hat{Y}_{B}}(\hat{f}_{B}(x),\hat{f}_{B}(y))\leq\hat{f}_{B}(x)\hat{f}_{B}(x')+\hat{f}_{B}(x')\hat{f}_{B}(y')+\hat{f}_{B}(y')\hat{f}_{B}(y).
$$

We already noted that the lengths of the segments outside $B$ are preserved by $\hat{f}_{B}$.
Therefore, it suffices to show that the length $\hat{f}_{B}(x')\hat{f}_{B}(y')$ is strictly smaller 
than $x'y'$ by a uniform amount. Thus, we can assume that the points $x$ and $y$ lie on distinct boundary components of $B$.

First, suppose that the images of the points $x$ and $y$ under $\hat{f}_{B}$ coincide,
that is, suppose that $x$ and $y$ are on the same equidistant curve from the core of $B$.
Then, 
$$
d_{\hat{X}}(x,y)-d_{\hat{Y}_{B}}(\hat{f}_{B}(x),\hat{f}_{B}(y))=d_{\hat{X}}(x,y)\geq\epsilon. 
$$
This proves the lemma in this case, and we are left to consider a pair of points $x$ and $y$ that do not lie 
on the same equidistant curve from the core of $B$, that is, curves whose images under $\hat{f}_{B}$
are  non-trivial geodesic segments.

Both projections of $[xy]$ on any side of $B$ along the arcs equidistant to the core give rise to a geodesic segment with 
exactly one endpoint in common with $[xy]$, and whose length equals that of $[\hat{f}_{B}(x)\hat{f}_{B}(y)$]. 
We now follow the same proof as that of Lemma \ref{lem:proj}, but taking care
this time of the difference between the length of the curve $[xy]$ and the length of its projection.
In order to lighten the reading, we keep the notation used in the proof of Lemma \ref{lem:proj}.
Thus, we denote the segment $[xy]$ by $\sigma$ and we shall specify a choice for the projection $k$ of 
$\sigma$.

Before explaining how to choose $k$, we remark that there is a  
upper bound $M=M(X,\alpha,\epsilon)>0$ to the length of such a projection $k$.
To see this, refer to Figure \ref{equid3}.
Any geodesic contained in $X$ and intersecting $B$ is contained in the
bi-infinite strip in $S_{\epsilon}$ bounded by $\widetilde{\gamma}_{1}$ and $\widetilde{\gamma}_{2}$.
This implies that there is an upper bound for the length of $\sigma$ that only depends
upon $X$ and on the choice of the $\epsilon$-strip $B$ in $\hat{X}$. 
Since the projection is $1$-Lipschitz, this implies the same result for the length of $k$.\\

If $\sigma$ intersects the core $c_{\epsilon}$ of $S_{\epsilon}$, then we subdivide $\sigma$ into 
two segments such that they both lie in different components of $S_{\epsilon}\setminus c_{\epsilon}$.
It then suffices to show the property for each of these segments.
We can therefore assume that $\sigma$ is contained in $S_{\epsilon}$ and that it does not intersect the core $c_{\epsilon}$ of $S_{\epsilon}$.
We are now back to the situation studied in the proof of Lemma \ref{lem:proj}, but with the
constraint on $\sigma$ of being contained in a strip of width $\epsilon$.
We keep the same choice for $k$ as the one settled in that proof and we 
refer the reader to Figure \ref{equid1} for what follows.

Consider the triangle $ABC$ as shown in Figure \ref{equid1}, but where  the three edges are now taken to be geodesics, namely, 
the edge $[AB]$ is the geodesic segment $\sigma$, the edge $[AC]$ is the geodesic segment $k$, and the third edge is the geodesic segment joining the points $B$ and $C$. Note that since the angle at $C$ of the triangle $\widehat{ACB'}$ of Figure 2 is equal to $\pi/2$, the angle $\widehat{ACB'}$ of the triangle that we consider now  is greater than or equal to $\pi/2$. We must show that the difference $AB-AC$ is bounded from
below by a uniform positive constant. 
Since the angle $\widehat{ACB}$ is greater or equal to $\pi/2$, we have
\begin{eqnarray*}
\cosh(AB)&=&\cosh(AC)\cosh(CB)-\sinh(AC)\sinh(CB)\cos(\widehat{ACB})\\
&\geq&\cosh(AC)\cosh(CB).
\end{eqnarray*}
Since $CB\geq\epsilon$, we get
$$
\cosh(AB)\geq\cosh(AC)\cosh(\epsilon).
$$
Now
\begin{eqnarray*}
\cosh(AB)&\geq&\cosh(AC)\cosh(\epsilon)\\
&\geq&\cosh(AC)(1+\epsilon^{2}/2).
\end{eqnarray*}
Hence,
\begin{eqnarray*}
\cosh(AB)-\cosh(AC)&\geq&\cosh(AC)\epsilon^{2}/2\\
&\geq&\epsilon^{2}/2.
\end{eqnarray*}
Multiplying the inequality by 2 and expanding $\cosh$ with exponentials, we get
\begin{eqnarray*}
e^{AB}-e^{AC}&\geq&\epsilon^{2}-(e^{-AB}-e^{-AC})\\
&\geq&\epsilon^{2}.
\end{eqnarray*}
The last inequality comes from the fact that $AB\geq AC$, that is, $e^{-AB}-e^{-AC}\leq0$.
We get
$$
e^{AB-AC}\geq 1+e^{-AC}\epsilon^{2},
$$
or,
$$
AB-AC\geq \log(1+e^{-AC}\epsilon^{2}).
$$
We saw that there exists a positive number $M=M(X,\alpha,\epsilon)>0$ such that $AC\leq M$.
Finally we get
$$
AB-AC\geq \log(1+e^{-M}\epsilon^{2})>0.
$$
This concludes the proof.
\end{proof}

The method used in the proof of  Proposition \ref{le:decrease} is due to Thurston \cite{Thurston1986}.

We need the following corollary in order to obtain the main result of this section (Theorem \ref{th:shrink} below).

We shall use the notion of measured geodesic lamination adapted to the case of surfaces with boundary, as in the paper \cite{Liu-Papad-Su-Th2009} and we briefly recall a few facts about geodesic laminations and their lengths. A measured geodesic lamination on the surface $S$ with boundary is defined in such a way that taking the double of this measured geodesic lamination  gives a measured geodesic lamination (in the usual sense) on the double of $S$, which is a surface without boundary. In particular, a measured geodesic lamination on $S$ is a finite union of uniquely defined minimal (with respect to inclusion) sub-laminations, called its {\it components}, and which are of the following three types:
\begin{enumerate}
\item A simple closed geodesic in $S$ (including a boundary component).
\item A geodesic arc meeting $\partial S$ at right angles.
\item A measured geodesic lamination in the interior of $S$, in which every leaf is dense. Such a component is called a {\it minimal component}.
\end{enumerate}

A measured geodesic lamination is said to be {\it  finite} if it has no minimal components.
The {\it length} of a measured geodesic lamination on $S$ is defined as in the case of surfaces without boundary, as the sum of the lengths of its components, and we recall a few facts about this notion. For a component which is an isolated leaf (that is, a simple closed geodesic or a geodesic arc), the length is the usual hyperbolic length. For a minimal component, one covers that component by geometric rectangles with disjoint interiors, and takes the sum of the areas of these rectangles. Here, a geometric rectangle (called, more simply, a {\it rectangle}) for the given measured geodesic lamination is  a quadrilateral immersed in $S$ having two opposite edges contained in leaves of the lamination and the other two edges transverse to the lamination. The length of such a rectangle is then defined as its total area, for the the area element defined as the product of the Lebesgue length element on the leaves of the lamination, and the one-dimensional measure on the transverse direction, provided by the transverse measure of the lamination. This is the notion of length of a lamination  used in Thurston's theory. A basic property of the length function is that it is continuous on the space of measured geodesic laminations on the surface. In particular, a property which is used several times  in \cite{Thurston1986} is that when a given measured lamination is approximated by a sequence of finite measured laminations, the lengths of these finite measured laminations converge to the length of the given measured lamination.

\begin{corollary}
Let $\lambda$ be a measured geodesic lamination on $S$.
Then $l_{Y_{B}}(\lambda)\leq l_{X}(\lambda)$, with strict inequality if and only if 
$\lambda\cap B\neq\emptyset$.
\end{corollary}

\begin{proof}
 First note that there is a natural correspondence between measured geodesic laminations on $X$ and measured geodesic laminations on $Y_B$, the two underlying surfaces being equivalent as marked surfaces. Thus, for any measured geodesic lamination $\lambda$ on $X$, we can talk about its length in $X$ and its length in $Y_B$. 
 
 Now if the support of $\lambda$ is a simple closed geodesic, the corollary follows from Proposition \ref{le:decrease}. If $\lambda$ is an arbitrary measured geodesic lamination, then, taking a sequence of weighted simple closed geodesics that approximate $\lambda$ in the topology of the space $\mathcal{ML}(S)$, the result follows from the continuity of the geodesic length function. In the case where $\lambda\cap B\not=\emptyset$, we can choose all the elements in the approximating sequence to satisfy this same property, and the result on the strict inequality also follows from Proposition \ref{le:decrease}. 
\end{proof}

\begin{theorem}\label{th:shrink}
For any point $X$ in Teichm\"uller space $\mathcal{T}(S)$, there exists a point $Y$ in $\mathcal{T}(S)$ such that $k(X,Y)<0$.
\end{theorem}

\begin{proof}
Choose a finite collection of geodesic arcs, $\mathcal{A}$, joining the boundary of $S$ to itself
such that any simple closed geodesic is intersected by one of these arcs.
Choose a collection of $\epsilon$-strip, one around each arc of $\mathcal{A}$.
Peel the  $\epsilon$-strips, one after the other.  
We thus get a new hyperbolic structure $Y$ on $S$ and a $1$-Lipschitz map
from the Nielsen extension of $X$ to that of $Y$.

Since any measured geodesic lamination of $X$ is intersected by an arc of $\mathcal{A}$,
the length of a measured geodesic lamination decreases when we pass from $X$ to $Y$.

Since $\mathcal{PML}(S)$ is compact, the supremum 
$$
\sup_{\alpha\in\mathcal{PML}(S)}\frac{l_{Y}(\alpha)}{l_{X}(\alpha)}
$$
is attained by a measured geodesic lamination.
Since the length of such a geodesic lamination has been strictly decreased, this shows that $k(X,Y)<0$.
This concludes the proof.
\end{proof}

\begin{remark}
The preceding result improves a theorem by Parlier \cite{Parlier} which says that for any
surface $S$ of finite type with non-empty boundary and for any hyperbolic structure $X$ on $S$, 
there exists a hyperbolic structure $Y$ on $S$ 
such that for every $\gamma$ in $\mathcal{C}$, we have $\displaystyle \frac{l_X(\gamma)}{l_Y(\gamma)}<1$. (Parlier's result only implies $k(X,Y)\leq0$). Note that whereas Parlier's result shows that the function $k$ defined in (\ref{eq:k}), in the case of a surface with nonempty boundary, is not a weak metric because it does not separate points, Theorem \ref{th:shrink} shows that for any surface with boundary, this function can even take negative values.
\end{remark}

\begin{remark}
Consider the peeling map $\hat{f}_{B}$ described above.
This map strictly decreases any elements of $\mathcal{B}\cup\mathcal{C}$ 
which intersects the strip $B$, and it leaves the lengths of 
the elements that are disjoint from $B$ unchanged. In the paper  \cite{Liu-Papad-Su-Th2009}, we defined the following function on 
the space  $\mathcal{T}(S)\times  \mathcal{T}(S)$ associated to a surface with boundary $S$:
$$d(X,Y)= \log \sup_{\gamma \in \mathcal{B}\cup\mathcal{C}}\frac{l_X(\gamma)}{l_Y(\gamma)}.
$$ 
and we showed that this function defines a weak metric. Thus, using Theorem \ref{th:shrink} above, there necessarily exists an arc on $S$ whose length increases when we pass from $X$ to $Y_{B}$. This arc is necessarily the arc $\alpha$ contained in the strip $B$.
It is therefore possible to compute the distance $d(X,Y_{B})$ explicitly.
\end{remark}

 We conclude with the following questions:
\begin{question}
Given a hyperbolic metric $X$ on a surface $S$ with nonempty boundary, can we always find another hyperbolic metric $Y$ such that every geodesic arc in $X$ which is length-minimizing between the boundary components is contracted when we pass from $X$ to $Y$ ? Note that the union of the arcs and boundary curves cannot all be contracted, by a result in \cite{Liu-Papad-Su-Th2009} that we already quoted above.
\end{question}

\begin{question}
We can ask the same question above, concerning geodesic arcs and interior geodesic closed curves, instead of only geodesic arcs.
\end{question}

\end{document}